\documentclass[a4paper,11pt,leqno]{amsart}
\usepackage[ansinew]{inputenc}
\usepackage{graphicx}
\usepackage{picins}

\usepackage{url}

\usepackage[colorlinks=true,linkcolor=blue,citecolor=magenta,urlcolor=red]{hyperref}

\hypersetup{pdfinfo={
  Title=On a Fair Copy of Riemann's 1859 Publication Created by Alfred Clebsch,
  Author=Wolfgang Gabcke,
  Subject=Handwritten copy of Riemann's 1859 Publication created by Alfred Clebsch,
  Keywords={Riemann's Nachlass, Bernhard Riemann, Alfred Clebsch, Richard Dedekind, Heinrich Weber,
            Karl Weierstraß, B. G. Teubner, Riemann zeta function, Riemann-Siegel formula}},
            linktocpage
}

\theoremstyle{plain}

\theoremstyle{definition}

\theoremstyle{remark}

\usepackage{scrextend}

\begin{document}
%
\title[ON A FAIR COPY OF RIEMANN'S 1859 PUBLICATION CREATED BY ALFRED CLEBSCH]
      {ON A FAIR COPY OF RIEMANN'S 1859\\PUBLICATION CREATED BY ALFRED CLEBSCH}
\author[WOLFGANG GABCKE]{WOLFGANG GABCKE}

\date{16th December 2015; 3rd release: Error messages (even linguistic) and suggestions are ap\-pre\-ci\-at\-ed}

\begin{abstract}
We will verify that the fair copy \textit{``Ueber die Anzahl der Primzahlen unter einer
gegebenen Grösse''} found in Riemann's Nachlass is not in Riemann's hand. Further, we will show that this paper
was written by Alfred Clebsch and that it was used after Clebsch's death with high probability as a setting copy
for chapter VII of Riemann's Collected Works published in 1876.
\end{abstract}
%
%
\maketitle
%
\deffootnotemark{\textsuperscript{(\thefootnotemark)}}
\deffootnote{2em}{1.6em}{(\thefootnotemark)\enskip}
\setcounter{tocdepth}{1}
\tableofcontents
%
%
%
%
%
%
\section{Introduction}
\noindent
Most of Bernhard Riemann's manuscripts, lecture notes, calculation sheets and letters that he left after his early
death in 1866 are preserved in 52 compilations at \textit{Nie\-der\-säch\-si\-sche Staats- und
Uni\-ver\-si\-täts\-bi\-blio\-thek Göt\-tin\-gen} (SUB) as ``Riemann's Nachlass''. One of these compilations is
one of the most famous compilations of mathematical original papers at all, it is Cod.\ Ms.\ B.\ Riemann 3
(Prim\-zah\-len), the Holy Grail of analytic number theory. Some of its papers are commonly
considered as the birth of this mathematical discipline. It contains
manuscripts connected with Riemann's 1859 publication \textit{``Über die An\-zahl der Prim\-zah\-len unter
einer ge\-ge\-be\-nen Größe''}\footnote{See \cite{Rie1, Rie2}; with modified title \textit{``Ueber die An\-zahl
der Prim\-zah\-len unter einer ge\-ge\-be\-nen Grösse''} see \cite[ch.\ VII, pp.\ 136--144]{Rie0} or \cite{Rie3}.
The printed version \cite{Rie1} was first published in spring 1860.}
as well as the derivation of the Riemann-Siegel integral formula, the Riemann-Siegel asymptotic formula, a lot of
calculation sheets in order to determine the first zeros of the zeta function, some drafts of letters to other
mathematicians\footnote{In a letter to Karl Weierstraß, Riemann mentions a new expansion of the zeta function
which is thought to be the Riemann-Siegel formula (see \cite{Rie6} and \cite[pp. 822--825]{Rie0a}).} and much
more. Apart from the 1859 publication, none of these findings were published by Riemann, although some of the
methods he used were new at that time, like the saddle point method in \cite{Rie7} to expand a complex integral in
an asymptotic series. The method remained unknown until 1909 when Peter Debye \cite{Deb} rediscovered it in
another unpublished paper of Riemann and used it for his aysmptotic expansions of the Bessel functions. The
integral formula and the asymptotic series mentioned above were unknown too, until C.\ L.\ Siegel
\cite{Sie1, Sie2} found them in compilation Cod.\ Ms.\ B.\ Riemann 3 in 1932, more than 70 years after Riemann had
derived them.

These opening words show that the Riemann hypothesis alone is not the reason, why mathematicians are
fascinated by Riemann's work about analytic number theory to this day. There are plenty of references
about the publication of 1859, the Riemann hypothesis and all other mathematical findings of compilation
Cod.\ Ms.\ B.\ Riemann 3 in the literature and on the Internet. Therefore, it is not surprising to find a
fac\-sim\-i\-le of a handwritten version of Riemann's 1859 paper\footnote{\label{Ref_Rie3}See \cite{Rie3} and figure
\ref{Abb03} for the first page.} from this compilation on the Internet, where it is thought that the original
manuscript of this fac\-sim\-i\-le not only was created but also written by Riemann himself as a fair copy. Although the
compilation came to SUB in 1890 and was inspected by many people since, no one noticed or was interested in the
fact that these six pages are not in Riemann's handwriting. This will be demonstrated in
section \ref{Copy_of_publication}. Further, we will show in section \ref{Alfred_Clebsch_creator} that the unknown
writer of the copy is Prof.\ Alfred Clebsch (1833--1872), successor to Riemann's Chair at the University of
Göttingen from 1868 to 1872. Then we will see in section \ref{Purpose_of_Copy} that the paper was
used as a setting copy for chapter VII of Riemann's Collected Works in 1876 with high probability.

To simplify the references, I have omitted ``Cod.\ Ms.\ B.\ Riemann'' in the following text; for example I will
write ``3 (19)'' instead of ``Cod.\ Ms.\ B.\ Riemann 3 (19)'' to reference the first page of Riemann's draft
\cite{Rie2}.

Some of the figures in the appendix are not suitable for printing on paper. They should be viewed on the screen
using the zoom function of a pdf reader.
%
%
%
%
%
%
\section{The Fair Copy of Riemann's Publication}
\label{Copy_of_publication}
\noindent
Manuscript 3 (16--18)\footnotemark[3] is a fair copy of Riemann's publication from 1859 not written in German
Kurrent\footnote{\label{Footnote_Kurrent}The Kurrent script was used in all German countries including some
regions in Switzerland, the Russian empire and in the Balkans as standard script for handwritten German text since
the 15th century. After a small simplification (Sütterlin) was introduced in 1915, 500 years of German script
tradition came to an end in 1941 when the Nazis banned its teaching in a circulare directive \cite{Rei} just like
the further use of its printed counterpart, the Fraktur. After the war this was not revised as
many other laws and decrees from the Nazi era too.} but in Latin script. Germans
usually wrote a German text in Kurrent script at that time but the creator of the fair copy did not. During the
examination of 14 compilations of Riemann's Nachlass I found no other person who wrote a German text in Latin
script. This applies for example to Ernst Abbe, Ri\-chard De\-de\-kind, Her\-mann Han\-kel, Karl Hat\-ten\-dorf,
Rie\-mann's widow Eli\-se, Her\-mann Aman\-dus Schwarz, Hein\-rich We\-ber, Karl Wei\-er\-straß and some more.
Indeed, the fair copy in Latin script seems to be unusual.

We will make a comparison of the handwritings between the fair copy 3 (16--18) and two autographs of Bernhard
Riemann. The first autograph is his draft 3 (19--20)\footnote{See \cite{Rie2} and figure
\ref{Abb01} for the first page.} of his communication to the \textit{Berliner Akademie} in October 1859. The date
can be seen at the end of line 2 on this page in the area with the water damage. The paper is written in Kurrent
script and the handwriting of its creator can be identified to be Riemann's by comparison with the large number of
preserved letters that contain Riemann's signature such as \cite{Rie5}. The second autograph is the fragment
3 (22)\footnote{See \cite{Rie4} and figure \ref{Abb02} for the first page.}, a fair copy of his publication in
French language and in Latin script, of which only one sheet (two pages) has survived.
Because of the different scripts, we cannot compare the German and the French text from 3 (19--20) and 3 (22)
directly to prove that fragment 3 (22) is in Riemann's hand. But we can compare mathematical formulas since their
notation was already internationally standardised back then, just as it is today. We give two examples in figure
\ref{Abb05} and \ref{Abb06} from page 3 (19) and 3 (22). It is not difficult to see that it must be the same
person who wrote these formulas as one compares the letters ``d'', ``n'', ``s'' and ``\(\Pi\)''. In both cases the
letter ``d'' looks like the notation of a partial derivative ``\(\partial\)'', a characteristic of Riemann's
handwriting. Also the style of letter ``s'' has a special feature in both papers: it often looks like digit ``5''.
Therefore, we can assume that the French fair copy 3 (22) must be in Riemann's hand.

Now, since 3 (22) is an autograph of Riemann in Latin script, we can compare the handwriting of the first page of
3 (16--18) in figure \ref{Abb03} with that of this autograph in figure \ref{Abb02}. We can immediately see that
the handwritings are completely different. Whereas the handwriting of 3 (16) must be that of a person who used it
every day since many years, Riemann's handwriting of 3 (22) shows an unpractised hand.
This is not surprising because he used Latin script only under special circumstances, namely when he wrote Latin,
French or English texts.

Comparing the style of some letters in figures \ref{Abb05} and \ref{Abb06}, the handwriting of mathematical
formulas do not match either. The creator of 3 (16) writes letter ``s'' as it is used today in contrast to
Riemann's ``5''. Further, his ``n'' and ``\(\Pi\)'' look different as well as his ``dx'' which does not contain a
``\(\partial\)''.

These comparisons show clearly that the handwriting of 3 (16) is not Riemann's so that 3 (16--18) cannot
have been written by Riemann himself.
%
%
%
%
%
%
\section{Alfred Clebsch, Creator of the Fair Copy}
\label{Alfred_Clebsch_creator}
\noindent
If Riemann is not the creator of the fair copy 3 (16--18), then who is? After a search through 13
other compilations of Riemann's Nachlass, I found a letter by Alfred Clebsch (1833--1872), successor to Riemann's
Chair at the University of Göttingen from 1868 to 1872, in compilation 1,2\footnote{See \cite{Cle1} or figure
\ref{Abb04} and subsection \ref{Transcription_of_letter} in the appendix for a German transcription and English
translation.}. Since this letter, written in Latin script and sent to the publishing company B.\ G.\ Teub\-ner at
Leipzig, is signed by Alfred Clebsch himself, it is an autograph and can be used for a comparison. If we do so and
look at figures \ref{Abb03} and \ref{Abb04}, we can immediately see that the handwritings of these documents are
identical. To substantiate this observation, it suffices to take a look at figure \ref{Abb07}--\ref{Abb12}, which
contains six capital letters from both documents. The result is clear: both documents must have been
written by the same person. Thus Alfred Clebsch is the creator of the fair copy 3 (16--18).
%
%
%
%
%
%
\section{The Presumed Purpose of the Fair Copy}
\label{Purpose_of_Copy}
\noindent
Alfred Clebsch was one of the editors of Riemann's Collected Works by Teubner at Leipzig until his death. The
other was Richard Dedekind. This can be gathered for example from Clebsch's letters \cite{Cle2} to Dedekind in
1872. In particular, we see it by his letter dated 17th of September 1872 \cite[no.\ 42]{Cle2} since it contains a
proposal for the designation of chapters and their order. Riemann's publication was provided under the title
``Primzahlen'' as chapter X. But after Clebsch's unexpected death on 7th November 1872, Heinrich Weber took over
the edition of Riemann's works and changed it to chapter VII\footnote{We cannot exclude that it was Clebsch
himself who made this change between 17th September and 25th October 1872 when he sent some of the setting copies
to Teubner (see his letter 2,5 (5) in \cite{Cle1} and in figure \ref{Abb04}). But this seems to be very unlikely.}.
This can be comprehended in the first line of Clebsch's fair copy 3 (16)\footnote{See \cite{Rie3} and \ref{Abb03}.}
where a Roman X was crossed and replaced by a Roman VII, all written with a pencil. Indeed, Riemann's publication
appeared as chapter VII of the Collected Works so Clebsch's fair copy must be connected to their edition in some
way.

Let us now look at the Collected Works \cite{Rie0} published in 1876 by Teubner at Leipzig, four years after
Clebsch's death. Since it is a scientific book, it was printed in Roman type and not in German
Fraktur\footnote{The ``Fraktur'' was the common print style for German text in German countries from the 15th
century until 1941. See footnote \ref{Footnote_Kurrent} on page \pageref{Footnote_Kurrent}.}. Although they are
defined in Roman type, the capital umlauts ``Ä'', ``Ö'' and ``Ü'' as well as the ``ß'' were not used in this book
but were replaced by ``Ae'', ``Oe'', ``Ue'' and ``ss''. This is exactly what Clebsch did in his letters to
Dedekind \cite{Cle2}, in his letter to Teubner \cite{Cle1}/figure \ref{Abb04} and in his fair copy 3 (16--18)
\cite{Rie3}/figure \ref{Abb03}, all written in Latin script. So was it Clebsch who had defined the print style of
the book? Definitively not, as can be seen in five other mathematical books \cite{Cle3, Plue, Sal, Schl, Ser}
published between 1868 and 1878 by Teubner that all have the same print style. So the commonly used title of
Riemann's publication
\begin{center}
\textit{\underline{Ue}ber die Anzahl der Primzahlen unter einer gegebenen Grö\underline{ss}e}
\end{center}
\noindent
is not correct since it is only due to Teubner's print style for mathematical books at that time. The
correct title
\begin{center}
\textit{\underline{Ü}ber die Anzahl der Primzahlen unter einer gegebenen Grö\underline{ß}e}
\end{center}
\noindent
should be used instead as Riemann defines it in his draft 3 (19--20) on the first line of figure \ref{Abb01} and as
it was printed in \textit{Monatsberichte der Akademie} \cite{Rie1}\footnote{The word ``über'' is used on page 671
instead of ``Über'' because the title is embedded in a larger text passage. But two pages before, on page 669, we
find the words ``Übersendung'' and ``Übersetzung'' showing that the print style of this book includes the capital
umlauts.} in 1860.

The different print style between the original publication \cite{Rie1} of 1860 and the reprint in the Collected
Works of 1876 \cite[ch.\ VII, pp.\ 136--144]{Rie0} gives reason to look at their deviations in detail, including
the handwritten versions of Riemann's draft 3 (19--20) \cite{Rie2} and Clebsch's fair copy 3 (16--18) \cite{Rie3}.
There are only minor deviations which are compiled in the table of figure \ref{Abb13}, with the limitation
that not all deviations between Riemann's draft and the print version of 1860 are listed. Up to three
non-significant clerical errors, there are only two deviations between Clebsch's fair copy and the print version of
the Collected Works concerning the words ``hiedurch''/``hierdurch'' and ``Hievon''/``Hiervon'' in deviation 1 and 3
of the table.
Clebsch first wrote ``Hiervon'' in deviation 3, then he corrected it to ``Hievon''. In the Collected Works we find
the linguistic modern spelling ``hierdurch'' and ``Hiervon'' which may be the result of a lector's revision. Thus
we can consider Clebsch's fair copy and the print version of 1876 as identical.

Deviation 4 is not significant since it is only a spelling variation, but deviation 2 is worthy of a more detailed
investigation. Riemann uses the word ``als''\ --\ English ``as''\ --\ in his draft which he also does in the print
version of 1860 in \textit{Monatsberichte der Akademie}. Clebsch, however, removes this word in his fair copy.
This change is also included in the Collected Works of 1876, four years after Clebsch's death. Since it is very
unlikely that another person, perhaps the new editor Heinrich Weber or a lector of the Teubner company, made
exactly the same modification of Riemann's text at the same location as Clebsch did four years earlier, there is
only one possible explanation: Clebsch's fair copy was used as a setting copy for chapter VII of Riemann's
Collected Works published in 1876, just as Er\-win Neu\-en\-schwan\-der suspected in a private communication.
Also, it is very likely that it belonged to the ``ordered compilation of prints and manuscripts'' that Clebsch
sent to Teubner at Leipzig on 25th October 1872\footnote{See subsection \ref{Transcription_of_letter} on page
\pageref{Transcription_of_letter}.}. Still, what we don't know is, how this document came into compilation
Cod.\ Ms.\ B.\ Riemann 3. Possibly, an answer can be found in the archives of the Teubner company, now preserved
in the \textit{Staatsarchiv Leipzig}.
%
%
%
%
%
%
\section*{Acknowledgements}
\noindent
I thank Bärbel Mund (SUB) and Nicolas Woodhouse for their suggestions, Helmut Rohlfing (SUB) for his comments and
text corrections, Erwin Neuenschwander for his hint to the crossed Roman X on Clebsch's fair copy and to Riemann's
letter to Weierstraß, Johannes Mangei (SUB) for his provision of advice and the motivation to write this essay and
Jonas Björklund as well as my children Wieland and Heidrun for their linguistic suggestions and corrections.
%
%
%
%
%
%

%
%
%
%
%
%
\newpage
\section{Appendix}
\subsection{Transcription of Clebsch's letter}
\label{Transcription_of_letter}
The letter is sheet 5 of compilation Cod.\ Ms.\ B.\ Riemann 1,2 (see \cite{Cle1} and figure \ref{Abb04}). Alfred
Clebsch wrote it on 25th October 1872, 13 days(!)\ before his tragic death of diphtheria on 7th November 1872, to
the publishing company B.\ G.\ Teubner at Leipzig in connection with the edition of Bern\-hard Rie\-mann's Collected
Works. This is indubitably an autograph of Alfred Clebsch since it is signed by himself. Presumably, this letter
is Clebsch's last preserved document at all.
\vspace{7ex}
\begin{center}
\large{\textsc{German Transcription}}
\end{center}
\vspace{1ex}
\noindent
\large{
\begin{tabbing}
\hspace{18.0ex}\=\hspace{1.7ex}\=\hspace{1.7ex}\=\hspace{3.9ex}\=
\hspace{5.6ex}\=\hspace{7.8ex}\=\hspace{6.7ex}\=\kill
\>\>\>\>\>\>Göttingen d.\ 25\raisebox{2mm}{\normalsize{\(\,.\,\)t}} Oct.\ 72.\\
\>An die Teubnersche Buchhandlung, Leipzig.\\\\
\>\>\>\>Hochgeehrter Herr!\\
\>\>Sie erhalten hiebei die zu der Ausgabe\\
\>von Riemanns Werken gehörigen Drucke\\
\>und Manuscripte geordnet; es fehlen nur\\
\>noch drei nicht sehr umfangreiche Ab-\\
\>handlungen, die ich ihrer Zeit folgen lasse.\\
\>\>\>Mit vollkommenster Hochachtung\\
\>\>\>\>\>\>Ihr\\\\
\>\>\>\>\>\>\>ergebenster\\
\>\>\>\>\>\>\>A.\ Clebsch
\end{tabbing}
}
\vspace{5ex}
\begin{center}
\large{\textsc{English Translation}}
\end{center}
\vspace{1ex}
\noindent
\large{
\begin{tabbing}
\hspace{18.0ex}\=\hspace{1.7ex}\=\hspace{1.7ex}\=\hspace{3.9ex}\=
\hspace{5.6ex}\=\hspace{7.8ex}\=\hspace{6.7ex}\=\kill
\>\>\>\>\>\>Göttingen, 25th Oct.\ 72.\\
\>To Teubner's bookstore, Leipzig.\\\\
\>\>\>\>Highly esteemed Sir!\\
\>\>Enclosed, you will find prints and manu-\\
\>scripts associated to the edition of Riemann's\\
\>works in an ordered compilation; there are\\
\>only three not very extensive treatises\\
\>missing that I let follow at their time.\\
\>\>\>With perfect tribute\\
\>\>\>\>\>\>your\\\\
\>\>\>\>\>\>\>obedient\\
\>\>\>\>\>\>\>A.\ Clebsch
\end{tabbing}
}
\newpage
\subsection{Figures}
\label{Figures}
\(\phantom{a}\)
\piccaption{\label{Abb01}Riemann's draft of his publication, first page, October 1859}
\parpic(\textwidth,117ex){\includegraphics[width=0.94\textwidth]{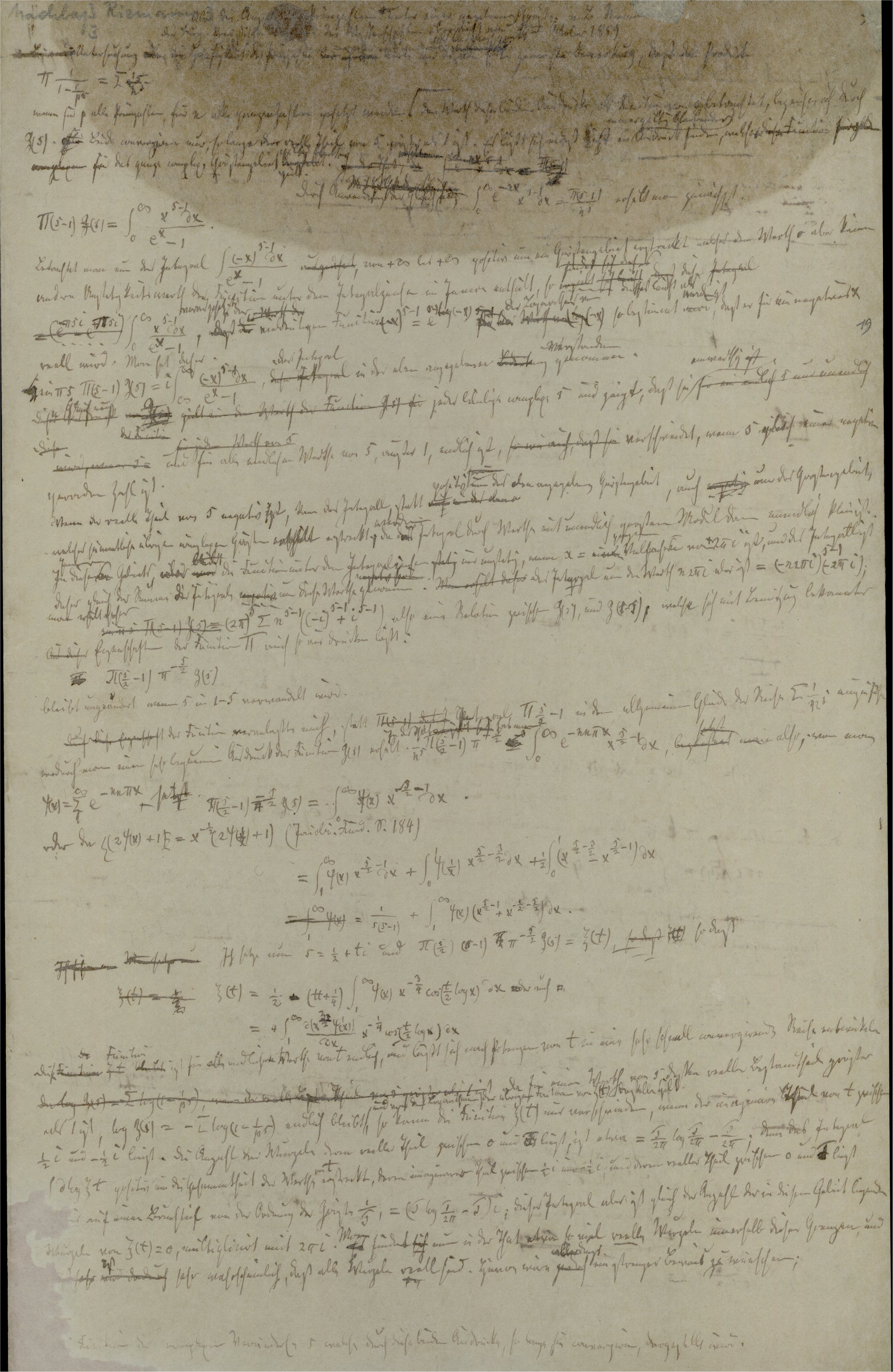}}
\newpage
\piccaption{\label{Abb02}Riemann's fair copy in French, first page}
\parpic(\textwidth,121ex){\includegraphics[width=0.90\textwidth]{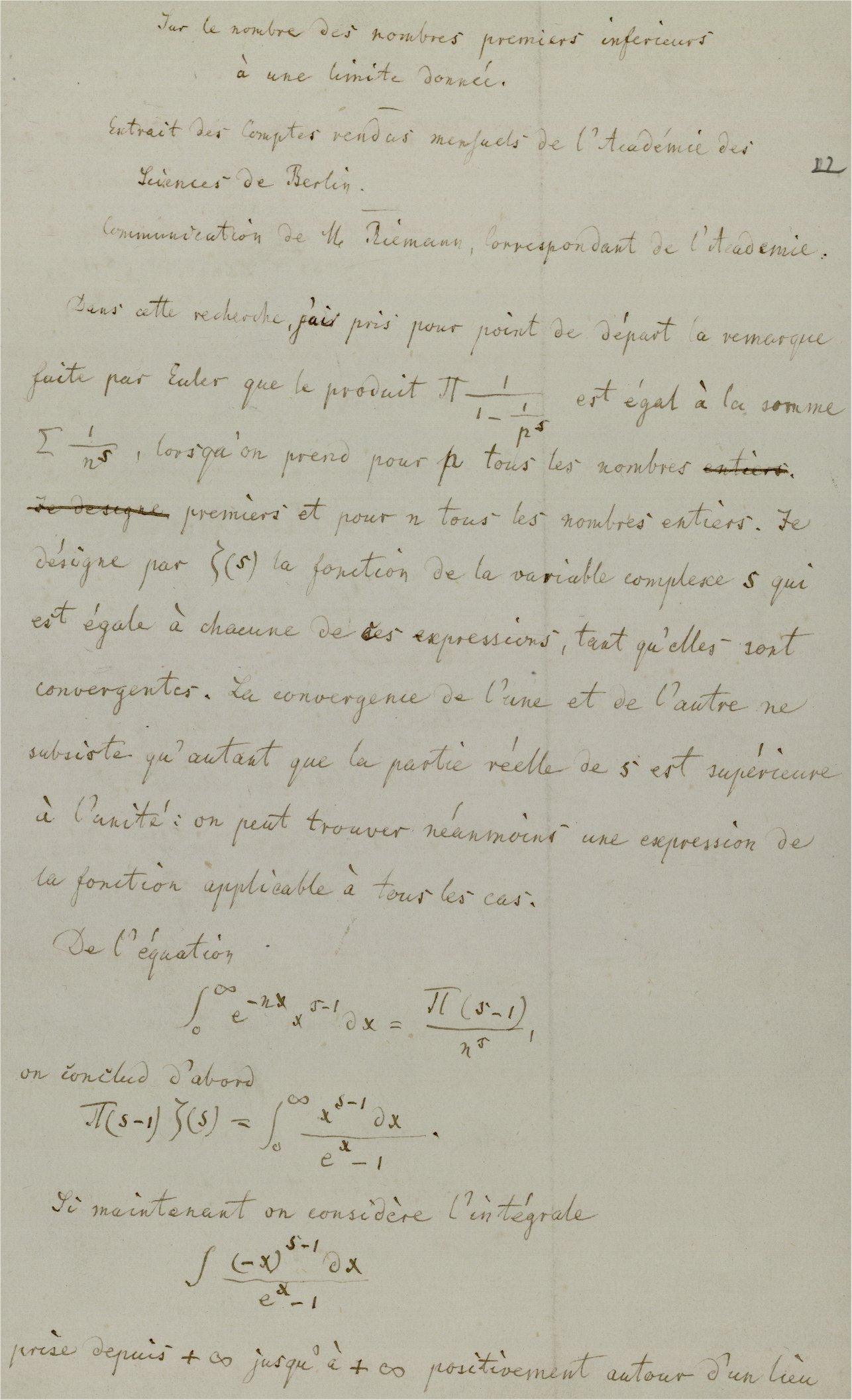}}
\newpage
\piccaption{\label{Abb03}Clebsch's fair copy of Riemann's publication, first page}
\parpic(\textwidth,121ex){\includegraphics[width=0.735\textwidth]{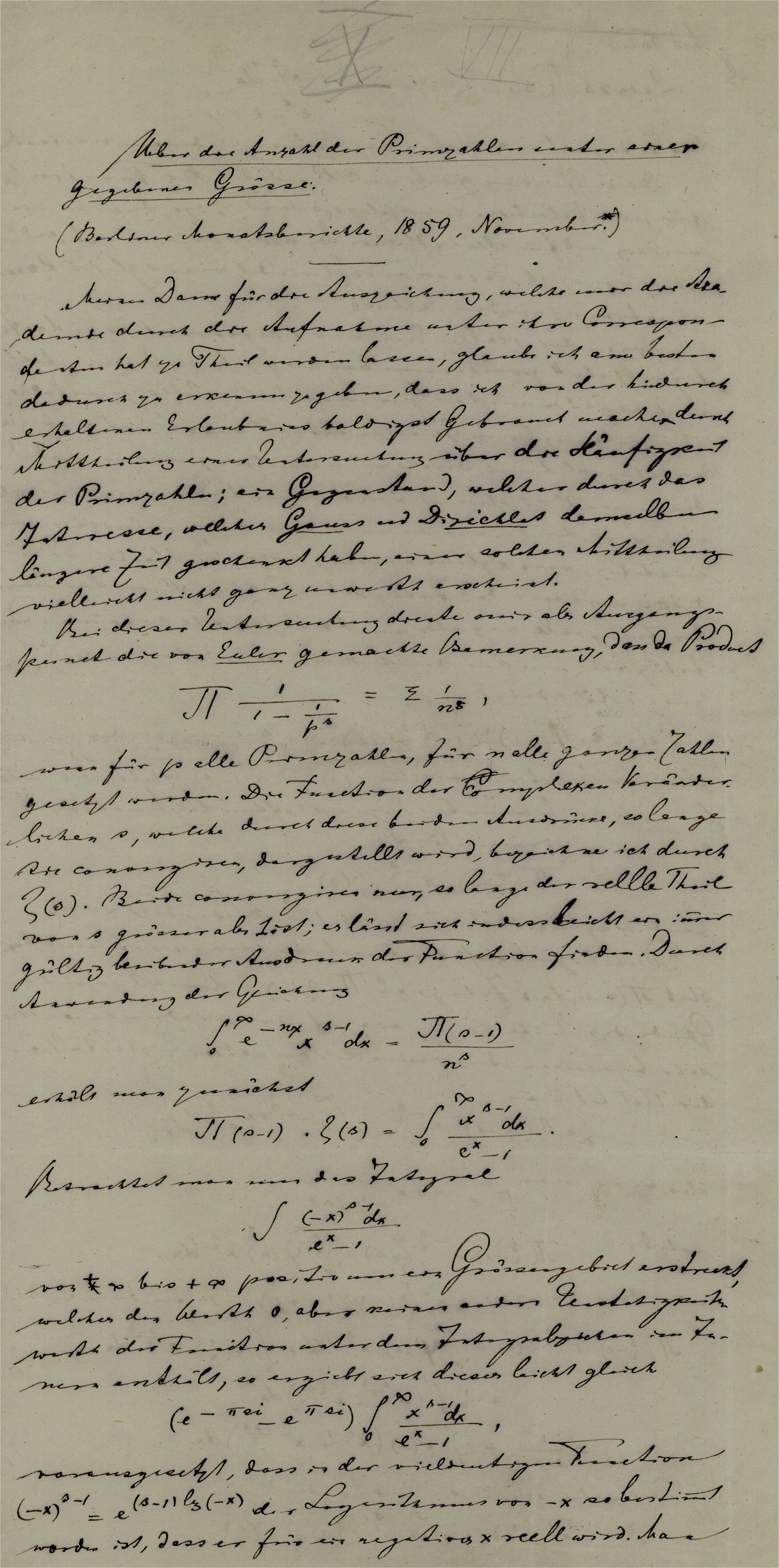}}
\newpage
\(\!\!\!\!\!\!\!\!\!\!\!\!\!\!\!\!\!\!\!\!\phantom{a}\)
\piccaption{\label{Abb04}Clebsch's letter to Teubner at Leipzig}
\parpic(\textwidth,95ex){\includegraphics[width=1.0\textwidth]{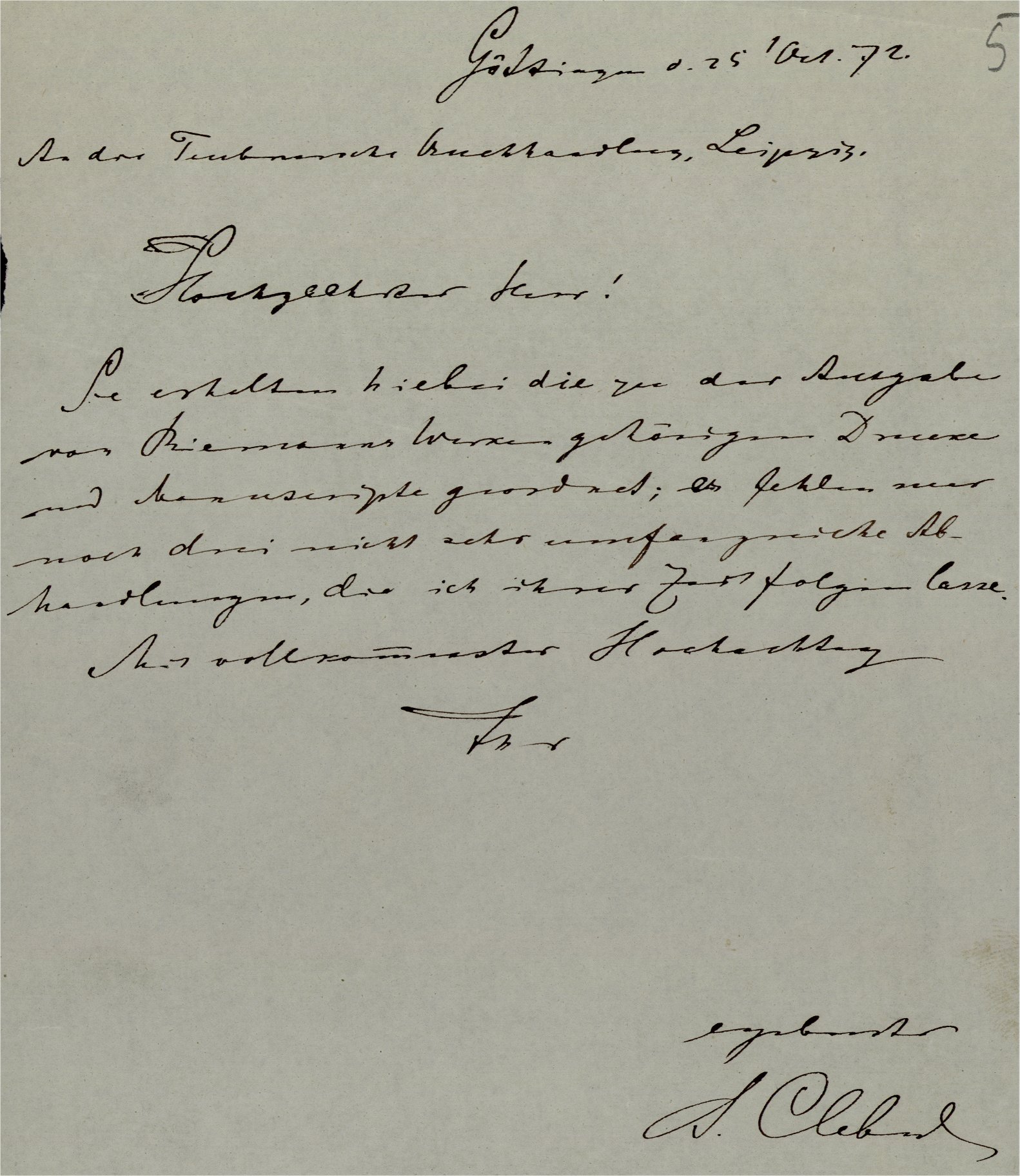}}
\newpage
\parpic(\textwidth,30.0ex){\includegraphics[width=0.85\textwidth]{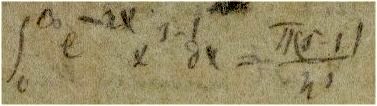}}
\vspace{20.5ex}
\parpic(\textwidth,21.7ex){\includegraphics[width=0.85\textwidth]{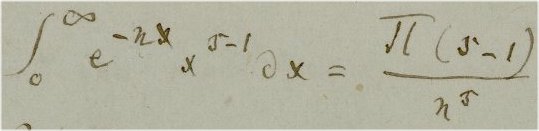}}
\vspace{15.2ex}
\piccaption{\label{Abb05}Details of 3 (19), 3 (22) and 3 (16) from top to bottom}
\parpic(\textwidth,17.7ex){\includegraphics[width=0.85\textwidth]{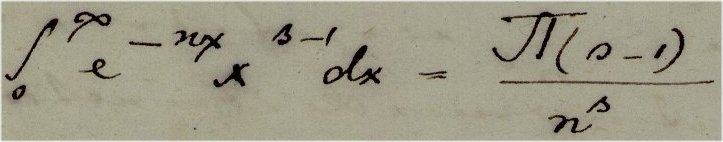}}
\vspace{27.6ex}
\parpic(\textwidth,19.0ex){\includegraphics[width=0.5\textwidth]{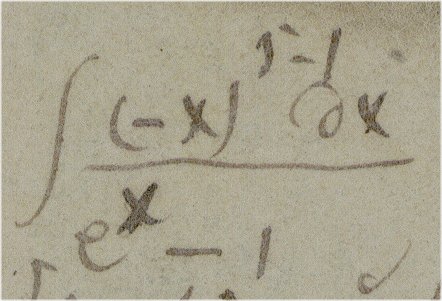}}
\vspace{19.9ex}
\parpic(\textwidth,21.6ex){\includegraphics[width=0.5\textwidth]{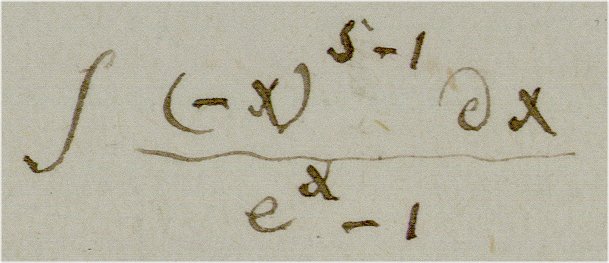}}
\vspace{16.2ex}
\piccaption{\label{Abb06}Details of 3 (19), 3 (22) and 3 (16) from top to bottom}
\parpic(\textwidth,20.1ex){\includegraphics[width=0.5\textwidth]{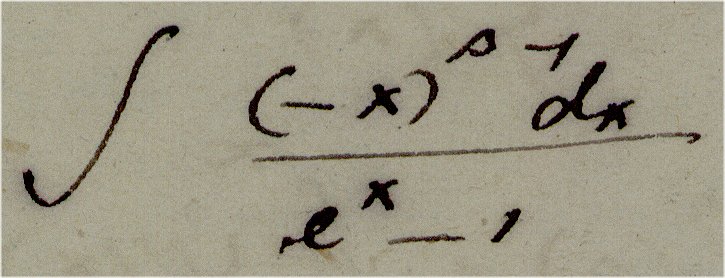}}
\newpage
\(\!\!\!\!\!\!\!\!\!\!\!\!\!\!\!\!\!\!\!\!\phantom{a}\)
\piccaption{\label{Abb07}``A'' in line 4 of 1,2 (5) and in line 5 of 3 (16)}
\parpic(72.0ex,20.5ex)(18.5ex,18.5ex){\includegraphics[width=0.25\textwidth]{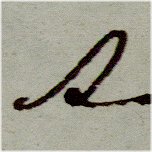}}
\parpic(36.0ex,20.5ex)(40.0ex,18.2ex){\includegraphics[width=0.25\textwidth]{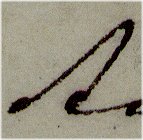}}
\vspace{25.0ex}
\piccaption{\label{Abb08}``D'' in line 5 of 1,2 (5) and in line 24 of 3 (16)}
\parpic(72.0ex,20.8ex)(18.5ex,18.5ex){\includegraphics[width=0.25\textwidth]{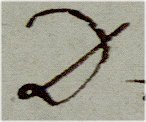}}
\parpic(36.0ex,20.8ex)(40.0ex,18.25ex){\includegraphics[width=0.20\textwidth]{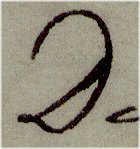}}
\vspace{46.5ex}
\piccaption{\label{Abb09}``G'' in line 1 of 1,2 (5) and in line 11/12 of 3 (16)}
\parpic(72.0ex,19.9ex)(18.5ex,17.6ex){\includegraphics[width=0.25\textwidth]{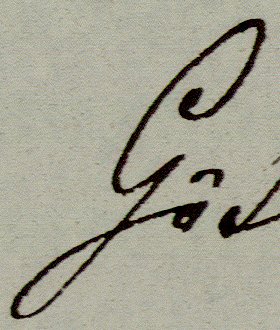}}
\parpic(36.0ex,19.9ex)(40.0ex,17.4ex){\includegraphics[width=0.25\textwidth]{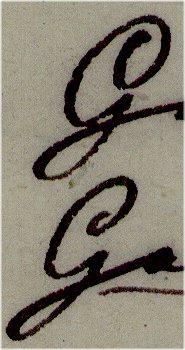}}
\vspace{24.5ex}
\piccaption{\label{Abb10}``L'' in line 2 of 1,2 (5) and in line 37 of 3 (16)}
\parpic(72.0ex,20.6ex)(18.5ex,18.5ex){\includegraphics[width=0.25\textwidth]{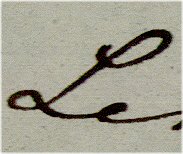}}
\parpic(36.0ex,20.6ex)(40.0ex,18.25ex){\includegraphics[width=0.28\textwidth]{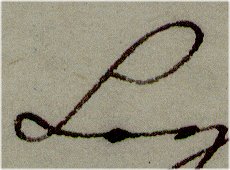}}
\newpage
\piccaption{\label{Abb11}``T'' in line 2 of 1,2 (5) and in line 7 of 3 (16)}
\parpic(72.0ex,20.6ex)(18.5ex,18.5ex){\includegraphics[width=0.25\textwidth]{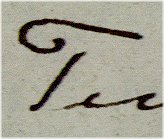}}
\parpic(36.0ex,20.6ex)(40.0ex,18.25ex){\includegraphics[width=0.30\textwidth]{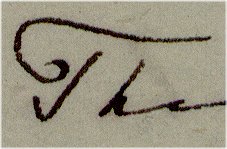}}
\vspace{28.5ex}
\piccaption{\label{Abb12}``Z'' in line 8 of 1,2 (5) and in line 13 of 3 (16)}
\parpic(72.0ex,20.8ex)(18.5ex,18.5ex){\includegraphics[width=0.25\textwidth]{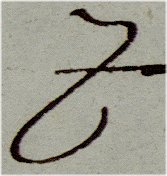}}
\parpic(36.0ex,20.8ex)(40.0ex,18.4ex){\includegraphics[width=0.20\textwidth]{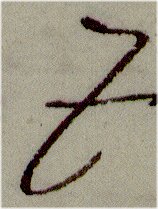}}
\vspace{23.0ex}
\piccaption{\label{Abb13}Deviations between the four versions of Riemann's publication}
\parpic(\textwidth,41.5ex){\includegraphics[width=1.25\textwidth]{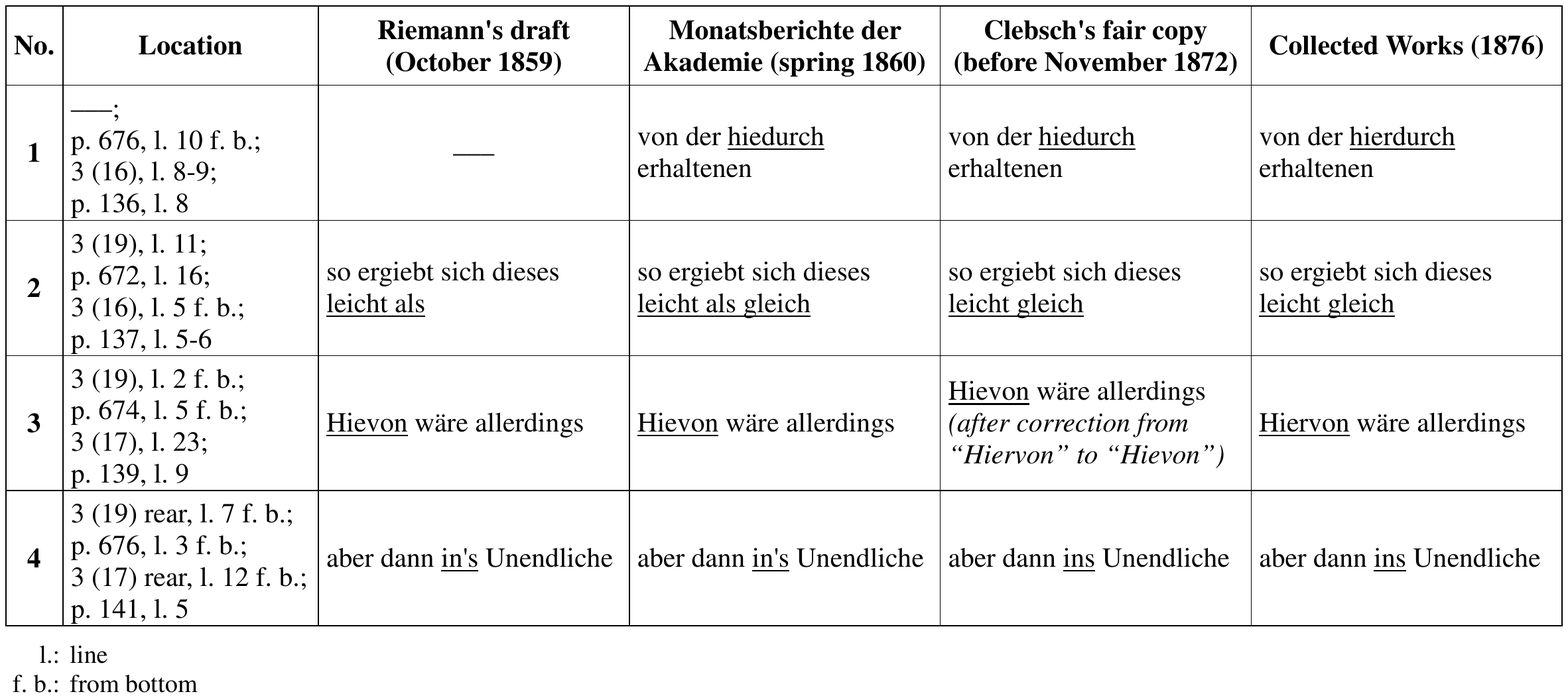}}
\end{document}